%&amstex          
\input amstex\documentstyle{amsppt}  
\pagewidth{12.5cm}\pageheight{19cm}\magnification\magstep1
\topmatter
\title Unipotent blocks and weighted affine Weyl groups\endtitle
\author G. Lusztig\endauthor
\address{Department of Mathematics, M.I.T., Cambridge, MA 02139}\endaddress
\thanks{Supported by NSF grant DMS-1855773.}\endthanks
\endtopmatter   
\document

\define\bW{\bar W}

\define\Irr{\text{\rm Irr}}

\define\si{\sim}

\define\sqc{\sqcup}

\define\part{\partial}
\define\emp{\emptyset}

\define\iy{\infty}
\define\m{\mapsto}
\define\do{\dots}

\define\lra{\leftrightarrow}

\define\sub{\subset}    

\define\T{\times}
\define\ti{\tilde}
\define\nl{\newline}
\redefine\i{^{-1}}

\define\un{\underline}

\define\ot{\otimes}

\define\Ad{\text{\rm Ad}}
\define\Hom{\text{\rm Hom}}

\define\Ind{\text{\rm Ind}}

\define\sg{\text{\rm sgn}}
\define\tr{\text{\rm tr}}

\redefine\b{\beta}
\redefine\c{\chi}
\define\g{\gamma}
\redefine\d{\delta}
\define\e{\epsilon}

\define\io{\iota}
\redefine\o{\omega}

\define\ph{\phi}

\define\s{\sigma}
\redefine\t{\tau}
\define\th{\theta}

\define\z{\zeta}
\define\x{\xi}

\redefine\G{\Gamma}
\redefine\D{\Delta}
\define\Om{\Omega}

\redefine\L{\Lambda}

\redefine\aa{\bold a}

\define\boc{\bold c}

\define\kk{\bold k}

\define\qq{\bold q}

\define\CC{\bold C}

\define\FF{\bold F}

\define\NN{\bold N}

\define\QQ{\bold Q}

\define\ZZ{\bold Z}

\define\ca{\Cal A}

\define\cc{\Cal C}

\define\ce{\Cal E}

\define\cg{\Cal G}

\define\cj{\Cal J}

\define\cl{\Cal L}

\define\cs{\Cal S}
\define\ct{\Cal T}

\define\cw{\Cal W}

\define\cx{\Cal X}

\define\fb{\frak b}
\define\fc{\frak c}

\define\fA{\frak A}

\define\fH{\frak H}

\define\fL{\frak L}

\define\fS{\frak S}
\define\fT{\frak T}

\define\fW{\frak W}

\define\tE{\ti E}

\define\sha{\sharp}

\define\bS{\bar S}

\head Introduction\endhead
\subhead 0.1\endsubhead
Let $G$ be a connected reductive algebraic group over $\CC$.
Let $un(G)$ be the set of unipotent conjugacy classes in $G$.
Let $ls(G)$ be the set of all pairs $(\fc,\fL)$ where $\fc\in un(G)$
and $\fL$ is an irreducible $\CC$-local system
on $\fc$, equivariant for the conjugation action of $G$.
In \cite{S76}, Springer discovered a remarkable bijection between a certain
subset ${}'ls(G)$ of $ls(G)$ and the set of irreducible representations of the
Weyl group of $G$ (up to isomorphism). In \cite{L84} I extended Springer result
by defining the ``generalized Springer correspondence'' that is,
a partition of $ls(G)$ into subsets (which could be called
{\it unipotent blocks}) and a bijection, for each unipotent block $\b$, between
the set of objects in $\b$ and the set $\Irr(\fW_\b)$
of irreducible representations (up to isomorphism) of a certain Weyl
group $\fW_\b$ associated to $\b$. (The subset ${}'ls(G)$ is one of the unipotent
blocks, namely the one containing $(\{1\},\CC)$.)
The arguments of \cite{L84} were based on the use of perverse sheaves on $G$.

\subhead 0.2\endsubhead
In the remainder of this introduction we assume that $G$ is almost simple, simply
connected. Let $W$ be the affine Weyl group of type dual to that of $G$.
In this paper we try to show that various information on unipotent elements of $G$
can be recovered from things which are more primitive
than unipotent elements, namely from Weyl groups, their associated 
Hecke algebras and their representations. We want to recover such information from
$W$ and its subgroups, without use of the geometry of $G$. Results of this type have been
obtained in \cite{L20} for the unipotent block ${}'ls(G)$ and this paper can be viewed as an
attempt to extend \cite{L20} to arbitrary unipotent blocks.

We now describe some earlier results in the same direction.

In \cite{L79} it has been observed that there is a set defined purely in terms
of $W$ which is an indexing set for $un(G)$, via the ordinary Springer
correspondence. More precisely, this set describes $un(G)$ in terms of truncated
induction from special representations of the various finite standard
parabolic subgroups of $W$, see 6.3. The same process allows one to recover
the order of the group of components of the centralizer of a unipotent element in
the adjoint group of $G$ in terms of data attached to the various special
representations in the previous sentence, see \cite{L09}.
In \cite{L80} it was conjectured (and in
\cite{L89} proved) that the set of two-sided cells of $W$ is an indexing set for
$un(G)$.
These results suggest that it may be possible to recover other properties
of unipotent elements of $G$ in terms of $W$.

\subhead 0.3\endsubhead
Let $\Om_W$ be the group of automorphisms of $W$ as a Coxeter group which
induce an inner automorphism of $W$ modulo the subgroup of translations.
In this paper, for any $\o\in\Om_W$, we define (without reference to $G$)
a certain set $\cc_\o(W)$ of standard finite parabolic subgroups $W_J$ of $W$,
stable under $\o$, see \S3.
Moreover, we define a bijection between the set of unipotent blocks of $G$ and the
set $\cup_{\o\in\Om_W}\cc_\o(W)$, see 6.5.
(This bijection is closely related to the arithmetic/geometric correspondence
[L95] in the study of unipotent representations of a simple $p$-adic group.)

One of the properties that we impose on $W_J\in\cc_\o(W)$ implies that a
finite reductive group over $\FF_q$ with Weyl group $W_J$ and with Frobenius
acting on this
Weyl group as $\o$ should have a unipotent cuspidal representation (but we
formulate this property without reference to finite reductive groups, see \S2).
This
has the consequence that the subgroups $W_J$ which appear are rather few.
In \S5 we attach to each $W_J\in\cc_\o(W)$ an affine Weyl group $\cw_J$
with a weight function $\cl_J:\cw_J@>>>\NN$. (This follows a procedure from \cite{L90},
\cite{L95}.) In this way to any unipotent block we have attached a weighted affine
Weyl group. (We regard the group $\{1\}$ as a weighted affine Weyl group.)
The weighted affine Weyl group $(\cw_J,\cl_J)$ is different in general from the one
defined in \cite{L17}.
In the remainder of this introduction
we  fix a unipotent block $\b$ and we denote by $\o,W_J,\cw_J,\cl_J$ the
objects associated to $\b$ as above.

In Theorem 6.10 we show that the quotient
$\bar\cw_J$ of $\cw_J$ by its group of translations
(which is defined purely in terms of $W$) can be identified with the group $\fW_\b$
of 0.1.
In \S4 we define a function $c:\Irr(\bar\cw_J)@>>>\NN$ purely in terms of $W$.
The values of $c$ are conjecturally related to dimensions of certain
Springer fibres (see 6.11); this relation is unconditional for exceptional types.
The function $c$ is used in 6.12(a) to express conjecturally the generalized
Green functions \cite{L86, 24.8} purely in terms of $W$; again this is
unconditional for exceptional types. In 7.3 we state a conjecture which
provides an indexing set for the set of two-sided cells in $\cw_J$ relative to
$\cl_J$; this is unconditional for exceptional types.

\subhead 0.4. Notation\endsubhead
For a group $L$ we denote by $L_{der}$ the derived subgroup of $L$.
For a finite group $\G$ let $\Irr(\G)$ be the set of irreducible representations of $\G$
over $\CC$ (up to isomorphism).

\head Contents\endhead
1. Weighted Weyl groups.

2. Sharp Weyl groups.

3. Affine Weyl groups and the sets $\cc_\o(W)$.

4. The function $c:\Irr(\bar\cw)@>>>\NN$.

5. A weighted affine Weyl group.

6. The set $\cg_\o(G)$ and the bijection $\cg_\o(G)@>\si>>\cc_\o(W)$.

7. Cells in the weighted affine Weyl group $\cw_J$.

\head 1. Weighted Weyl groups\endhead
\subhead 1.1\endsubhead
Let $\fW$ be a Coxeter group and let $\fS$ be its set of simple reflections.
Le $w\m|w|$ be the length function on $\fW$.
We say that $\fW$ is {\it weighted} if we are given a weight function $\cl:\fW@>>>\NN$
that is a function such that $\cl(ww')=\cl(w)+\cl(w')$ for any $w,w'$ in $\fW$ such that
$|ww'|=|w|+|w'|$.
Let $Cell(\fW,\cl)$ be the set of two sided cells of $(\fW,\cl)$ in the sense
of \cite{L03}.

Let $v\in\CC-\{0\}$ be a non-root of $1$ and let $\fH_v$ be the Hecke algebra
(over $\CC$) associated to $\fW,\cl,v$; thus $\fH_v$ is the $\CC$-vector space with basis
$\{T_w;w\in\fW\}$ with associative multiplication defined by the rules
$T_wT_{w'}=T_{ww'}$ for
any $w,w'$ in $\fW$ such that $|ww'|=|w|+|w'|$ and $(T_\s+v^{-\cl(\s)})(T_\s-v^{\cl(\s)})=0$
for any $\s\in\fS$.

\subhead 1.2\endsubhead
In the remainder of this section we assume that $\fW$ is a Weyl group with
a given weight function $\cl:\fW@>>>\NN$.
Let $\Irr(\fH_v)$ be the set of simple $\fH_v$-modules
(up to isomorphism). It is known that $\Irr(\fW),\Irr(\fH_v)$ are in canonical
bijection.
(We can assume that $\fW$ is irreducible. If properties P1-P15 in \cite{L03, \S14} are assumed,
then according to \cite{L03}, the algebras
$\fH_v,\CC[\fW]$ are canonically isomorphic and the result would follow. Now P1-P15 do hold
when $\fW$ is of type $E_7,E_8$ or $G_2$. Thus we can assume that $\fW$ is of type other than
$E_7,E_8$ or $G_2$. In this case, the result 
follows from the observation in \cite{BC} according to which the various $E\in\Irr(\fW)$ are
characterized by their multiplicities in representations of $\fW$ induced from the unit
representation of various parabolic subgroups. Alternatively we can use \cite{G11}.) We denote this
bijection by $E\lra E_v$. For $E\in\Irr(\fW)$ we set
$$f_{E,\cl,v}=(\dim E)\i\sum_{w\in\fW}\tr(T_w,E_v)\tr(T_{w\i},E_v)\in\CC^*$$
and
$$D_{E,\cl,v}=f_{E,\cl,v}\i\sum_{w\in\fW}v^{2\cl(w)}\in\CC.$$
From the known explicit formulas for $D_{E,\cl,v}$ we see that $D_{E,\cl,v}$
is a nonzero rational function in $v$. We define $\aa_\cl(E)\in\NN$ by the requirement that
$D_{E,\cl,v}v^{-2\aa_\cl(E)}$ is a rational function in $v$ whose value at $v=0$ is $\ne0$
and $\ne\iy$. 
In the case where $\cl=||$ we write $\aa(E)$ instead of $\aa_\cl(E)$.

\subhead 1.3\endsubhead
If $\fW$ is of type $A_1$ and the value of $\cl_{\fS}$ is $a\in\ZZ_{>0}$ then setting $q=v^{2a}$ we
have $D_{E,\cl,v}=1$ if
$E=1$ and $D_{E,\cl}=q$ if $E=\sg$; the corresponding values of $\aa_{E,\cl}$ are $0;a$.
We now assume that $\fW$ is of type $A_2$ and the values of $\cl|_{\fS}$ are $a,a$
in $\ZZ_{>0}$. Setting $q=v^{2a}$ we write below the values of $D_{E,\cl,v}$ for various $E$
of dimension $1;2;1$:
$$1;q^2+q,q^3.$$
The corresponding values of $\aa_{E,\cl}$ are $0;a;3a$.

\subhead 1.4\endsubhead
In this subsection we assume that $\fW$ is of type $B_2$ and the values of $\cl|_{\fS}$ are $a,b$
in $\ZZ_{>0}$. Setting $q=v^{2a},y=v^{2b}$ we write below the values of $D_{E,\cl,v}$ for various $E$
of dimension $1;2;1;1;1$:
$$1;qy(q+1)(y+1)/(q+y);q^2(qy+1)/(q+y);y^2(qy+1)/(q+y);q^2y^2.$$
The corresponding values of $\aa_{E,\cl}$ are:
$$0;a+b-m;2a-m;2b-m;2a+2b$$
where $m=\min(a,b)$.

\subhead 1.5\endsubhead
In this subsection we assume that $\fW$ is of type $G_2$ and the values of $\cl|_{\fS}$ are $a,b$
in $\ZZ_{>0}$. Setting $q=v^{2a},y=v^{2b}, \sqrt{qy}=v^{a+b}$ we write below the values of
$D_{E,\cl,v}$ for various $E$ of dimension $1;2;2;1;1;1$:
$$\align&
1;qy(q+1)(y+1)(qy+\sqrt{qy}+1)/(2(q+\sqrt{qy}+y));\\&
qy(q+1)(y+1)(qy-\sqrt{qy}+1)/(2(q-\sqrt{qy}+y));\\&
q^2(q^2y^2+qy+1)/(q^2+qy+y^2);y^2(q^2y^2+qy+1)/(q^2+qy+y^2);\\&
q^3y^3.\endalign$$
The corresponding values of $\aa_{E,\cl}$ are:
$$0;a+b-m;a+b-m;2a-2m;2b-2m,3a+3b$$
where $m=\min(a,b)$.

\subhead 1.6\endsubhead
In this subsection we assume that $\fW$ is of type $B_3$ and the values of $\cl|_{\fS}$ are $a,a,b$
in $\ZZ_{>0}$. Setting $q=v^{2a},y=v^{2b}$ we write below the values of $D_{E,\cl,v}$ for various $E$
of dimension $1;3;2;1;3;3;2;3;1;1$:
$$\align&1;qy(q^2+q+1)(qy+1)/(q+y);q^2(q+1)(q^2y+1)/(q+y);\\&
y^3(qy+1)(q^2y+1)/((q^2+y)(q+y));\\&
qy^2(q^2+q+1)(q^2y+1)/(q^2+y);q^3y(q^2+q+1)(q^2y+1)/(q^2+y);\\&
q^2y^3(q+1)(q^2y+1)/(q+y);q^3y^2(q^2+q+1)(qy+1)/(q+y);\\&
q^6(qy+1)(q^2y+1)/((q^2+y)(q+y));q^6y^3.\endalign$$

The corresponding values of $\aa_{E,\cl}$ are:
$$\align&0;a+b-m;2a-m;3b-m-m';a+2b-m';3a+b-m';2a+3b-m;\\&
3a+2b-m;6a-m-m';6a+3b\endalign$$
where $m=\min(a,b)$, $m'=\min(2a,b)$.

\subhead 1.7\endsubhead
In \cite{L82} a definition of a partition of $\Irr(\fW)$ into subsets called
families was given. Repeating that definition but using $\aa_\cl(E)$ instead of
$\aa(E)$ for $E\in\Irr(\fW)$ we obtain a partition of $\Irr(\fW)$ into subsets called
$\cl$-families. (This definition appears in \cite{L83, no.7}.) Thus the families
of \cite{L82} are the same as the $||$-families.

\head 2. Sharp Weyl groups\endhead
\subhead 2.1\endsubhead
Let $\fW,\fS,||$ be as in 1.1. We assume that $\fW$ is a Weyl group.
Let $A_{\fW}$ be the group of all automorphisms $\g$ of $\fW$ preserving $\fS$ and such that
whenever $\s\ne\s'$ are in the same $\g$-orbit in $\fS$, the product $\s\s'$ has 
order $\ge3$. For $\g\in A_{\fW}$ let $r(\g)$ be the number of $\g$-orbits on
$S$; let $ord(\g)$ be the order of $\g$. We shall also write ${}^\g\fW$ instead
of $(\fW,\g)$. (We sometimes
write ${}^d\fW$ instead of ${}^\g W$ where $d=ord(\g)$;
when $d=1$ we write $\fW$ instead of ${}^1\fW$.)

Let $op\in A_{\fW}$ be given by conjugation by the longest element of $\fW$. 

For any $E\in\Irr(\fW)$ and $v\in\CC^*$ a non-root of $1$, $D_{E,||,v}\in\CC$ is defined as in 1.2 with $\cl=||$.
It is known that $D_{E,||,v}$ is a polynomial in $v^2$ with rational coefficients.
Let $z(E)$ be the largest integer $\ge0$ such that $D_{E,||,v}/(v^2+1)^{z(E)}$ is a
polynomial in $v^2$. From the known formulas for $D_{E,||,v}$ one can see that $z(E)\le r(op)$.

Assuming that $\fW$ is $\{1\}$ or irreducible, we say that $\fW$ is {\it sharp} if

(a) there exists $E_0\in\Irr(\fW)$ such that $z(E_0)=r(op)$
\nl
and $r(op)ord(op)$ is even. (The last condition is imposed to rule out a $\fW$ of type $E_7$.)
Note that the family of $\fW$ containing $E_0$ in (a) is necessarily unique (when
such $E_0$ exists); moreover $E_0$ itself is unique (when it exists) if it is assumed to be special.

Let $\g\in A_{\fW}$. Assuming that $\fW$ is $\{1\}$ or irreducible, we say that
${}^\g\fW$ is sharp if $\fW$ is sharp and $ord(op\g)$ is odd. The last condition
means that we have either $\g=op$ or $ord(\g)=3$ (hence $\fW$ is of type $D_4$). 
For ${}^\g\fW$  sharp we set $\aa[{}^\g\fW]=\aa(E_0)$ where $E_0$ is as in (a).

\subhead 2.2\endsubhead
Here is a complete list of the various sharp ${}^\g\fW$ and the corresponding
$\aa[{}^\g\fW]$:

(i) $\{1\}$, $\aa[{}^\g\fW]=0$:

(ii) ${}^2A_{(t^2-1)/8-1}$, $t\in\{5,7,9,\do\}$,
$\aa[{}^\g\fW]=(t-3)(t-1)(t+1)/48$;

(iii) $B_{(t^2-1)/4}$, $t\in\{3,5,7,\do\}$, $\aa[{}^\g\fW]=(t-1)(t+1)(2t-3)/24$;

(iv) $D_{t^2/4}$, $t\in\{4,8,12,\do\}$, $\aa[{}^\g\fW]=(t-2)t(2t+1)/24$;

(v) ${}^2D_{t^2/4}$, $t\in\{6,10,14,\do\}$, $\aa[{}^\g\fW]=(t-2)t(2t+1)/24$;

(vi) $G_2$, ${}^3D_4$, $F_4$, ${}^2E_6$, $E_8$,
$\aa[{}^\g\fW]=1,3,4,7,16$ respectively.
\nl
For a general $\fW,\g$ we say that $\fW$ is $\g$-irreducible if $\fW$ is a product of
$k\ge1$ irreducible Weyl groups $\fW_1,\do,\fW_k$ and $\g$ permutes $\fW_1,\do,\fW_k$
cyclically. In this case  we say that ${}^\g\fW$ is sharp if ${}^{\g^k}\fW_1$ is
sharp and we set $\aa[{}^\g\fW]=k\aa[{}^{\g^k}\fW_1]$.

 For a general $\fW,\g$, we have that  $\fW$ is a product $\fW'_1\T\do\T\fW'_l$ of Weyl groups
such that each $\fW'_j$ is $\g$-stable and $\g$-irreducible. We say that ${}^\g\fW$
is sharp if each ${}^\g\fW'_j$ is sharp; we set
$\aa[{}^\g\fW]=\sum_j\aa[{}^\g\fW_j]$.

The objects in (i),(iii)-(v) can be viewed as vertices of a graph:
$$\align&\{1\} ---B_{(3^2-1)/4}--- D_{4^2/4}---B_{(5^2-1)/4}---\\&
---{}^2D_{6^2/4}---B_{(7^2-1)/4}---D_{8^2/4}---B_{(9^2-1)/4}---\\&
---{}^2D_{10^2/4}---B_{(11^2-1)/4}---
\do.\tag a\endalign$$
   
We will attach to each vertex of this graph an index: the index of each of
$B_{(t^2-1)/4}$, $D_{t^2/4}$, ${}^2D_{t^2/4}$ is $t$; the index of $\{1\}$ is $2$.

From the objects in (i)-(v) we can form a second graph:
$$\align &\{1\}\T\{1\}---{}^2A_{(5^2-1)/8-1}---
{}^\g(B_{(6/2)^2-1)/4}\T B_{((6/2)^2-1)/4})---\\&
---{}^2A_{(7^2-1)/8-1}---{}^\g(D_{(8/2)^2/4}\T D_{(8/2)^2/4})---\\&
---{}^2A_{(9^2-1)/8-1}---
{}^\g(B_{((10/2)^2-1)/4}\T B_{((10/2)^2-1)/4})---\\&
---{}^2A_{(11^2-1)/8-1}---{}^\g(D_{(12/2)^2/4}\T D_{(12/2)^2/4})---\\&
---{}^2A_{(13^2-1)/8-1}---\do. \tag b \endalign$$

Here $\g$ acts on $B_{((t/2)^2-1)/4}\T B_{((t/2)^2-1)/4}$ as an involution
exchanging the two factors; it acts on 
$D_{(t/2)^2/4}\T D_{(t/2)^2/4}$ by permuting the two factors in such a way
that $\g^2=1$ if $t\in\{8,16,24,\do\}$ and $\g^2\ne1$ if $t\in\{12,20,28,\do\}$.
We will attach to each vertex of this graph an index: the index of each of 
$${}^\g(D_{(t/2)^2/4}\T D_{(t/2)^2/4}),{}^2A_{(t^2-1)/8-1},
{}^\g(B_{(t/2)^2-1)/4}\T B_{((t/2)^2-1)/4})$$
is $t$; the index of $\{1\}\T\{1\}$ is $4$.

\head 3. Affine Weyl groups and the sets $\cc_\o(W)$\endhead
\subhead 3.1\endsubhead
In this section $W$ denotes an (irreducible) affine Weyl group.
Let $\ct$ be the set of all $w\in W$ such that the conjugacy class of $w$ is
finite. (Such $w$ are said to be the translations of $W$.)
Now $\ct$ is a free abelian group of finite rank and of finite index in $W$.
Let $w\m|w|$ be the usual length function of $W$. Let $S$ be
the set of simple reflections of $W$. Let $S^!$ be the set of all $\s\in S$ such
that the sum of labels of edges of the Coxeter graph of $W$ which touch $\s$ is
$\ge3$. We have $\sha(S^!)\le2$. 

Let $\Om_W$ be the (finite abelian) group of automorphisms of $W$ preserving $S$
whose restriction to $\ct$ is given by conjugation by an element of $W$. 
If $S^!\ne\emp$ let $\Om'_W$ be the set of all $\o\in\Om_W$ such that
$\o$ restricted to $S^!$ is the identity map (this is a subgroup of $\Om_W$);
if $S^!=\emp$ we set $\Om'_W=\Om_W$. We set $\Om''_W=\Om_W-\Om'_W$.

Let $\bW=W/\ct$ (a finite group). We show:

(a) {\it If $\o\in\Om_W$ then $\o:W@>>>W$ induces an inner automorphism of $\bW$.}
\nl
We can find $w\in W$ such that $\o(\t)=\Ad(w)(\t)$ for all $\t\in\ct$.
Let $\z=\Ad(w\i)\o:W@>>>W$. We have $\z(\t)=\t$ for any $\t\in\ct$.
Let $y\in W,\t\in\ct$. We have $y\t y\i\in\ct$ hence $\z(y\t y\i)=y\t y\i$ that is
$\z(y)\t\z(y\i)=y\t y\i$. Setting $y'=y\i\z(y)\in W$ we have
$y'\t y'{}\i=\t$ for any $\t\in W$. Now the action of $W/\ct$ on $\ct$ by
conjugation is faithful hence $y'\in\ct$. Thus $\z(y)\in y\ct$ so that
$w\i\o(y)w\in y\ct$ and $\o(y)\in wyw\i\ct$. This proves (a).

\subhead 3.2\endsubhead
For any $J\subsetneqq S$ let $W_J$ be the subgroup of $W$ generated by $J$ (a finite Weyl group).

Let $S_*$ be the set of all $\s\in S$ such that $W_{S-\{\s\}}@>>>\bW$
(restriction of the obvious map $W@>>>\bW$) is an isomorphism. We have $S_*\ne\emp$
and the obvious action of $\Om_W$ on $S_*$ is simply transitive. If $J\subsetneqq S$ let
$W_J@>>>\bW$ be the restriction of the obvious homomorphism $W@>>>\bW$; this is an
imbedding, so that $W_J$ can be viewed as a subgroup of $\bW$.
For any special representation $E\in\Irr(W_J)$
there is a unique $E'\in\Irr(\bW)$ such that
$E'$ appears in $\Ind_{W_J}^{\bW}(E)$ and in the $\aa(E)$-th symmetric power of
the conjugation representation of $\bW$ on $\CC\ot\ct$ (with $\aa(E)$ defined in
terms of $W_J$); we set $E'=j_{W_J}^{\bW}(E)$, see \cite{L09, 1.3}.
Let $\bS(\bW)$ be the subset of $\Irr(\bW)$ consisting of representations of
the form $j_{W_J}^{\bW}(E)$ for some $J\subsetneqq S$ and some special 
$E\in\Irr(W_J)$.

\subhead 3.3\endsubhead
Let $\o\in\Om_W$. We define a set $\cc_\o(W)$ of Weyl subgroups $W_J$ of $W$
with  $J\subsetneqq S$. This set contains $\{1\}$. Now $\cc_\o(W)-\{1\}$ consists
of the subgroups $W_J$ with $J\subsetneqq S$, $J\ne\emp$
which satisfy the following requirements.

(i) $W_J$ is $\o'$-stable for any $\o'\in\Om_W$.

(ii) $W_J$ is $\o$-sharp.

(iii) $W_{S-J}$ is $\o$-irreducible.

(iv) If $\sha(S^!)=2$ and $\Om''_W=\emp$ then  $W_J$ is the product of two
vertices of the graph 2.2(a) which are joined by an edge.

(v) if $\o\in\Om''_W$, then $W_J$ is the product of two vertices of the graph
2.2(b) which are joined by an edge.

\subhead 3.4\endsubhead
We now describe the set $\cc_\o(W)$ in each case.
If $W$ is of affine type $A_{n-1}$, $n\ge2$, then $\Om_W$ is cyclic of order $n$;
for any $\o\in\Om_W$, $\cc_\o(W)$ consists of a single element: $\{1\}$
(with $\aa[\{1\}]=0$).

\subhead 3.5\endsubhead
If $W$ is of affine type $E_6$ then $\Om_W$ is cyclic of order $3$.
Let $\o\in\Om_W$. If $\o=1$ then $\cc_\o(W)$ consists of a single element: $\{1\}$
(with $\aa[\{1\}]=0$).
If $\o\ne1$ then $\cc_\o(W)$ consists of $\{1\}$ (with $\aa[\{1\}]=0$) and
of the subgroup $W_J$ of
type $D_4$ (so that ${}^\o W_J={}^3D_4$ and $\aa[{}^\o W_J]=3$).

\subhead 3.6\endsubhead
If $W$ is of affine type $E_7$ then $\Om_W$ is cyclic of order $2$.
Let $\o\in\Om_W$. If $\o=1$ then $\cc_\o(W)$ consists of a single element: $\{1\}$
(with $\aa[\{1\}]=0$).
If $\o\ne1$ then $\cc_\o(W)$ consists of $\{1\}$ (with $\aa[\{1\}]=0$)
and of the subgroup $W_J$ of
type $E_6$ (so that ${}^\o W_J={}^2E_6$ and $\aa[{}^\o W_J]=7$).

\subhead 3.7\endsubhead
If $W$ is of affine type $E_8,F_4$ or $G_2$ and $\o\in\Om_W$ then $\o=1$
and $\cc_\o(W)$ consists of $\{1\}$ (with $\aa[\{1\}]=0$)
and of the subgroup
$W_J$ of non-affine type $E_8,F_4$ or $G_2$ (respectively), with $\aa[W_J]$ equal
to $16,4,1$ respectively.

\subhead 3.8\endsubhead
In the remainder of this section we assume that $W$ is of affine type
$B_n (n\ge3),C_n (n\ge2)$ or $D_n (n\ge4)$.
Let $\o\in\Om_W$. Assume first that $\o\in\Om'_W$.
Let ${}'\cc_\o(W)$ be the set of all pairs $(t,s)\in\NN^2$ such that

$t-s=\pm1$ (type $B$), $t=s$ (type $C,D$),

$t=0\mod4$ (type $B,D$ with $\o=1$), $t=2\mod4$ (type $B,D$ with $\o\ne1$)
$t=1\mod2$ (type $C$),

and for some $r\in\NN$ we have

$t^2/4+(s^2-1)/4+r=n$, (type $B$), 

$(t^2-1)/4+(s^2-1)/4+r=n$, (type $C$),

$t^2/4+s^2/4+r=n$, (type $D$).

that is, $ts+2r=2n$    (type $B,D$),  $ts+2r=2n+1$ (type $C$).

In type $B$ we define a bijection $\cc_\o(W)@>\si>>{}'\cc_\o(W)$ by associating to
$W_J\in\cc_\o(W)$ (assumed to be $\ne\{1\}$) the pair $(t,s)$ formed by the
indexes $t,s$ of the two vertices
attached to $W_J$ in 3.3(iv) and by associating to $W_J=\{1\}$ the pair $(0,1)$
(if $\o=1$) or $(2,1)$ (if $\o\ne1$).

In type $C,D$, we define a bijection $\cc_\o(W)@>\si>>{}'\cc_\o(W)$ by
associating to $W_J\in\cc_\o(W)$ (assumed to be $\ne\{1\}$) the pair $(t,s)$ (with
$t=s$) where $W_J$ is the product of a vertex of index $t$ in 2.2(a) with itself
and by associating to $W_J=\{1\}$ the pair $(1,1)$ (type $C$), $(0,0)$ (type $D$
with $\o=1$), $(2,2)$ (type $D$ with $\o\ne1$).

Assuming that $W_J$ corresponds as above to $(t,s)$ we can compute
$\aa[{}^\o W_J]$ in each case.

If $s=t\in\{1,3,5,\do\}$ then 
$$\align&\aa[{}^\o W_J]=(t-1)(t+1)(2t-3)/24+(t-1)(t+1)(2t-3)/24\\&
=(t-1)(t+1)(2t-3)/12.\endalign$$
If $s=t\in\{0,2,4,\do\}$ then 
$$\aa[{}^\o W_J]=(t-2)t(2t+1)/24+(t-2)t(2t+1)/24=(t-2)t(2t+1)/12.$$
If $t\in\{0,2,4,\do\},s\in\{1,3,5,\do\},t-s=\pm1$ then
$$\aa[{}^\o W_J]=(t-2)t(2t+1)/24+(s-1)s(2s+1)/24$$
and this equals $(t-1)t(t+1)/6$ if $s=t+1$ and $(t-2)(t-1)t/6$ if $t=s+1$.

\subhead 3.9\endsubhead
Assume next that $\o\in\Om''_W$. Let ${}'\cc_\o(W)$ be the set of all pairs
$(t,s)\in\NN^2$ such that

$t-s=\pm1$, 

$t=2\mod4$ (type $C$),

$t=0\mod8$ (type $D$ with $\o^2=1$), $t=4\mod8$ (type $D$ with $\o^2\ne1$),

and for some $r\in\NN$ we have

$((t/2)^2-1)/2+(s^2-1)/8+2r=n$ (type $C$),

$(t/2)^2/2+(s^2-1)/8+2r=n$ (type $D$),

or equivalently $ts/2+4r=2n+1$ (type $C$), $ts/2+4r=2n$ (type $D$).

Note that $(2,1)\in{}'\cc_\o(W)$ in type $C$ with $n$ even and
$(2,3)\in{}'\cc_\o(W)$ in type $C$ with $n$ odd.
We define a bijection $\cc_\o(W)@>\si>>{}'\cc_\o(W)$ by associating to
$W_J\in\cc_\o(W)$ (assumed to be $\ne\{1\}$)
the pair $(t,s)$ formed by the indexes of the two vertices attached to
$W_J$ in 3.3(v) and by associating to $W_J=\{1\}$ the pair $(2,1)$
(type $C$ with $n$
even), the pair $(2,3)$ (type $C$ with $n$ odd), the pair $(0,1)$ (type $D$ with
$\o^2=1$), the pair $(4,3)$ (type $D$ with $\o^2\ne1$).

Assuming that $W_J$ corresponds as above to $(t,s)$ we can compute
$\aa[{}^\o W_J]$ in each case.

If $t\in\{2,6,10,\do\}$ then 
$$\align&\aa[{}^\o W_J]=(s-3)(s-1)(s+1)/48+(t/2-1)(t/2+1)(t-3)/24+\\&
(t/2-1)(t/2+1)(t-3)/24\endalign$$
which equals $(t-2)(2t^2-5t-6)/48$ if $s=t-1$ and equals $(t-2)(t+2)(2t-3)/48$ if
$s=t+1$.

If $t\in\{0,4,8,10,\do\}$ then 
$$\aa[{}^\o W_J]=(s-3)(s-1)(s+1)/48+(t/2-2)(t/2)(t+1)/24+(t/2-2)(t/2)(t+1)/24$$
which equals $(t-4)t(2t-1)/48$ if $s=t-1$ and equals $t(2t^2-3t-8)/48$ if $s=t+1$.

\subhead 3.10\endsubhead
Let $\o\in\Om'_W$. Let ${}'\un\cc_\o(W)$ be the set of all pairs $(\d,r)\in\NN^2$ such that

(type $B$)  $\d+2r=2n$, $\d=2+4+6+\do+(2\s)$, where $\s\in\NN$ and
$\s=0\mod4$ or $\s=3\mod4$ if $\o=1$; $\s=1\mod4$ or $\s=2\mod4$ if $\o\ne1$;

(type $C$) $\d+2r=2n+1$, $\d=1+3+5+\do+(2\s-1)$ where $\s\in\NN$ and $\s=1\mod2$;

(type $D$) $\d+2r=2n$, $\d=1+3+5+\do+(2\s-1)$, where $\s\in\NN$ and
$\s=0\mod4$ if $\o=1$, $\s=2\mod4$ if $\o\ne1$.

We define a bijection
$${}'\cc_\o(W)@>\si>>{}'\un\cc_\o(W)$$
by $(t,s)\m(ts,r)$ where $r\in\NN$ is as in 3.8.

\subhead 3.11\endsubhead
Let  $\D=\{t(t+1)/2;t\in\NN\}$. Let $\o\in\Om''_W$. 
Let ${}'\un\cc_\o(W)$ be the set of all pairs $(\d,r)\in\D\T\NN$ such that

$\d+4r=2n+1$ (type $C$), 

$\d+4r=2n$, $\d=0\mod4$ (type $D,\o^2=1$),

$\d+4r=2n$, $\d=2\mod4$ (type $D,\o^2\ne1$).
\nl
We define a bijection
$${}'\cc_\o(W)@>\si>>{}'\un\cc_\o(W)$$
by $(t,s)\m(ts/2,r)$ where $r\in\NN$ is as in 3.9.

\head 4. The function $c:\Irr(\bar\cw)@>>>\NN$\endhead
\subhead 4.1\endsubhead
In this section $\cw$ denotes an irreducible affine Weyl group with a set $\cs$
of simple
reflections and with a given weight function $\cl:\cw@>>>\NN$.
Let $||$ be the length function of $\cw$.
Let $\ct_\cw$
be the group of translations of $\cw$ (see 3.1) and let $\bar\cw=\cw/\ct_\cw$
(a finite group). For any $\cj\subsetneqq\cs$ we denote by $\cw_\cj$ the subgroup
of $\cw$ generated by $\cj$ (a finite Weyl group).
Let $\bar\cw_\cj$ be the image of $\cw_\cj$ under the obvious map $\cw@>>>\bar\cw$.
Note that the obvious map $\cw_\cj@>>>\bar\cw_\cj$ is an
isomorphism; we use this to identify $\cw_\cj=\bar\cw_\cj$.

For any $E\in\Irr(\bar\cw)$ we define $\fT(E)$ to be the set of all pairs
$(\cj,E')$ where $\cj\sub\cs$, $\sha(\cj)=\sha(\cs)-1$ and
$E'\in\Irr(\cw_\cj)$ is such that $E'$
appears in the restriction of $E$ to $\cw_j=\bar\cw_\cj$. We set
$$c_E=\max_{(\cj,E')\in\fT(E)}\aa_{\cl}(E')\in\NN,\tag a$$
$$\fT^*(E)=\{(\cj,E')\in\fT(E);\aa_{\cl}(E')=c_E\},$$
where $\aa_{\cl}(E')$ is defined as in 1.2 in terms of the Weyl group $\cw_\cj$
with the weight function obtained by restricting $\cl:\cw@>>>\NN$ to 
$\cw_\cj$. We have $\fT^*(E)\ne\emp$.

\subhead 4.2\endsubhead
The function $E\m c_E$ has been computed explicitly in \cite{L20}
in the case where $\cw$
is of exceptional type and $\cl=||$. We will now describe explicitly the map
$\Irr(\bar\cw)@>>>\NN,E\m c_E$ of 4.1 in several examples with $\cl\ne||$.
Assume that either

(a) $\cw$ is of affine type $C_r,r\ge2$ and $\cl|_{\cs}$ takes the values
$t,1,1,\do,1,1,s$ where $t>0,s>0$, or

(b) $\cw$ is of affine type $C_r,r\ge2$ and $\cl|_{\cs}$ takes the values
$t,2,2,\do,2,2,s$ where $t>0,s>0$, or

(c) $\cw$ is of affine type $G_2$ and $\cl|_{\cs}$ has values $3,3,1$.

In case (a),(b) we set $u=\max(t,s)$.
In each of the examples below we give a table with two rows whose columns are
indexed by the various $E\in\Irr(\bar\cw)$. The first row represents the numbers
$c_E$. The second row represents the numbers $\aa_{\cl}(E')$ for various
$E'\in\Irr(\cw_{\cj})$ where $\cj=\cs-\{s\}$ for some $s\in\cs_*$ (see
3.2) which
in case (a),(b) is chosen so that $\cl(s)=u$. (We can identify
$\cw_\cj=\bar\cw$ hence $E'$ can be identified with an $E\in\bar\cw_J$.)
Any entry $e$ in the first row and column $E$ is $\ge$ than the 
entry $e'$ in the second row and column $E$. When $e>e'$ we indicate some other
$\cj'\subsetneqq\cs$ such that some 
$E'\in\Irr(\cw_{\cj'})$ appears in the restriction of $E$ to $\cw_{\cj'}$ and
$\aa_{\cl}(E')=e$. (We will specify such $\cj'$ by specifying the type
of $\cw_{\cj'}$.)

Assume first (in case (a)) that $r=2$ and $t=s=u\ge1$. The table is
$$\align&0;1;u;2u;2u+2\\&
         0;1;t;2u-1;2u+2  \endalign$$
with an additional $\cj'$ with $\cw_{\cj'}$ of type $A_1\T A_1$,
contributing $2u$ to the fourth column.

Assume next (in case (a)) that $r=2$ and $t=s\pm1,u\ge2$. The table is
$$\align&0;1;u;2u-1;2u+2\\&
         0;1;u;2u-1;2u+2  \endalign$$
In this case there is no need for an additional $\cj'$.

Assume next (in case (b)) that $r=2$ and $t=s\pm1$, $u\ge2$. The table is
$$\align&0;u;2;2u-1;2u+2\\&
         0;u;2;2u-2;2u+2  \endalign$$
with an additional $\cj'$ with $\cw_{\cj'}$ of type $A_1\T A_1$,
contributing $2u-1$ to the fourth column.

Assume next (in case (b)) that $r=3$ and $t=s\pm1$, $u\ge2$. The table is
$$\align&0;u;2;3u-3;2u-1;u+2;3u+3;2u+4;6;3u+12\\&
         0;u;2;3u-6;2u-2;u+2,3u+2;2u+4;6;3u+12\endalign$$
with additional $\cj'$ with $\cw_{\cj'}$ of type $A_1\T B_2$,
contributing $3u-3$ to the fourth column, $2u-1$ to the fifth column and
$3u+3$ to the seventh column.

Finally assume that we are in case (c). The table is
$$\align& 0;1;3;4;9;12\\&
          0;1;3;3;7;12\endalign$$
with an additional $\cj'$ with $\cw_{\cj'}$ of type $A_2$
contributing $9$ to the fifth column and with an additional $\cj'$
with $\cw_{\cj'}$ of type $A_1\T A_1$ contributing $4$ to the fourth column. 

\subhead 4.3\endsubhead
We set
$$\nu(\cw,\cl)=\max_{s\in \cs_*}\cl(w_{0,s})$$
where $w_{0,s}$ is the element of maximal length of $\cw_{\cs-\{s\}}$.

\subhead 4.4\endsubhead
Let $\si$ be the equivalence relation on $\Irr(\bar\cw)$ generated by the
relation $E_1\si E_2$ when $c_{E_1}=c_{E_2}$ and there exist 
$(\cj_1,E'_1)\in\fT^*(E_1)$, $(\cj_2,E'_2)\in\fT^*(E_2)$ such that
$\cj_1=\cj_2$ and $E'_1,E'_2$ are in the same $\cl$-family (see 1.7) of
$\Irr(\cw_{\cj_1})=\Irr(\cw_{\cj_2})$.

\head 5. A weighted affine Weyl group\endhead
\subhead 5.1\endsubhead
We preserve the setup of 3.1.
For any $J\subsetneqq S$ let $w_0^J$ be the longest element
of $W_J$. We now fix $\o\in\Om_W$, $J\subsetneqq S$ such that $W_J\in\cc_\o(W)$.
Following \cite{L90} (see also \cite{L95}), to $W,\o,J$ we associate
a weighted affine Weyl group $(\cw_J,\cl_J)$. (A similar procedure was used for
finite Weyl groups in \cite{L78}.)
The proofs of various statements in this section can be extracted from \cite{L95}.
Let
$$\cw'_J=\{w\in W;wW_J=W_Jw,w\text{ has minimal length in }wW_J=W_Jw\};$$
this is a subgroup of $W$ stable under $\o$; let $\cw_J=\{w\in\cw'_J;\o(w)=w\}$ 
(another subgroup of $W$). Let $\cs_J$ be the set of all $\o$-orbits on $S-J$.

When $\sha(\cs_J)\ge2$ let $\ct_J$ be the group of translations of the affine Weyl group $\cw_J$ and let $\bar\cw_J=\cw_J/\ct_J$, a finite group.
When $\sha(\cs_J)=1$ we have $\cw_J=\{1\}$ and we set
$\ct_J=\{1\},\bar\cw_J=\{1\}$; in this case we regard $\cw_J$ as a weighted
affine Weyl group with weight function $\cl_J=0$.

We now assume that $\sha(\cs_J)\ge2$. If $\th\in\cs_J$ we set
$\t_\th=w_0^{J\cup \th}w_0^J=w_0^Jw_0^{J\cup \th}$.
(The last equality is a property of any $W_J\in\cc_\o(W)$.)
We have $\t_\th\in\cw_J,\t_\th^2=1$.
Moreover $\cw_J$ is a Coxeter group on the generators
$\{\t_\th;\th\in\cs_J\}$ and with Coxeter relations $(t_\th t_{\th'})^{m_{\th,\th'}}=1$
for any distinct $\th\in\cs_J,\th'\in\cs_J$ such that $J\cup\th\cup\th'\ne S$,
where
$$m_{\th,\th'}=\frac{2(|w_0^{J\cup\th\cup\th'}|-|w_0^J|)}
{|w_0^{J\cup\th}|+|w_0^{J\cup\th'}|-2|w_0^J|}.$$

(When $J\cup\th\cup\th'=S$ then $\t_\th\t_{\th'}$ has infinite order.)
Note that $\cw_J$ is an (irreducible) affine Weyl group.

\subhead 5.2\endsubhead
We define a weight function $\cl_J:\cw_J@>>>\NN$. Let $\th\in\cs_J$.
Let $E_0$ be the unique irreducible special representation of $W_J$ such that
$z(E_0)=r(op)$ (notation of 2.1 with $\fW$ replaced by $W_J$). Let
$$\ce=\{\tE\in\Irr(W_{J\cup\th};\tE\text{ appears in }\Ind_{W_J}^{W_{J\cup\th}}(E_0)\},$$
$$\cl(\th)=\max\{\aa(\tE);\tE\in\ce\}-\min\{\aa(\tE);\tE\in\ce\}.$$
(Here $\aa(\tE)$ is as in 1.2 with $\fW$ replaced by $W_{J\cup\th}$.)
This defines the function $\cl_J:\cs_J@>>>\ZZ_{>0}$. This extends uniquely to a
weight function $\cw_J@>>>\NN$ denoted again by $\cl_J$ so that $(\cw_J,\cl_J)$ is
a weighted affine Weyl group.

\subhead 5.3\endsubhead
Let $\o\in\Om_W$, $J\subsetneqq S$ be such that $W_J\in\cc_\o(W)$.
We now describe the pair $(\cw,\cl)=(\cw_J,\cl_J)$ in each case.
We will write $\cs$ instead of $\cs_J$.

If $W$ is of affine type $A_{n-1},n\ge2$, let $k=ord(\o)$ (a divisor of $n$).
We have $W_J=\{1\}$. If $k<n$, $\cw$ is of affine type $A_{(n/k)-1}$ with
$\cl|_{\cs}$ being constant equal to $k$. If $k=n$, we have $\cw=\{1\}$. 

Assume that $W$ is of affine type $E_6$. If $\o=1$ then $W_J=\{1\}$ and
$\cw=W$ with $\cl|_{\cs}$ being constant equal to $1$. If $\o\ne1$ and $W_J=\{1\}$ then
$\cw$ is of affine type $G_2$ with the values of $\cl|_{\cs}$ being $3,3,1$.
If $\o\ne1$ and $W_J$ is of type $D_4$ then $\cw=\{1\}$.

Assume that $W$ is of affine type $E_7$. If $\o=1$ then $W_J=\{1\}$ and
$\cw=W$ with $\cl|_{\cs}$ being constant equal to $1$. If $\o\ne1$ and $W_J=\{1\}$ then
$\cw$ is of affine type $F_4$ with the values of $\cl|_{\cs}$ being $2,2,2,1,1$.
If $\o\ne1$ and $W_J$ is of type $E_6$ then $\cw=\{1\}$.

Assume that $W$ is of affine type $E_8,F_4$ or $G_2$. We have $\o=1$.
If $W_J=\{1\}$ then $\cw=W$ with $\cl|_{\cs}$ being constant equal to $1$. If
$W_J$ is of the non-affine type $E_8,F_4$ or $G_2$ (respectively) then $\cw=\{1\}$.

\subhead 5.4\endsubhead
We now assume that $W$ is of affine type $B_n (n\ge3),C_n (n\ge2)$ or $D_n (n\ge4)$.
Let $(\d,r)\in{}'\un\cc_\o(W)$ be the image of $W_J$ under the composition
$\cc_\o(W)@>\si>>{}'\cc_\o(W)@>\si>>{}'\un\cc_\o(W)$
(the first map as in 3.8, 3.9, the second map as in 3.10, 3.11).

If $r=0$ then $\cw=\{1\}$.
If $r>0,\d>0$, then $\cw$ is of affine type $C_r$, with Coxeter graph
$$\boxed{t}==\boxed{1}--\boxed{1}--\do--\boxed{1}--\boxed{1}==\boxed{s} \tag a$$
$$\boxed{t}==\boxed{2}--\boxed{2}--\do--\boxed{2}--\boxed{2}==\boxed{s} \tag b$$
where the boxed entries are the values of $\cl|_{\cs}$ and $t,s$ in $\NN$ are
defined by

$ts=\d$, $t-s\in\{0,1,-1\}$ (in (a) with $\o\in\Om'_W$),

$ts/2=\d$, $t-s\in\{1,-1\}$ (in (b) with $\o\in\Om''_W$).
\nl
(Here affine of type $C_1$ is taken to be the same as affine of type $A_1$.)

If $\d=0$ (hence $r>0$), then either:

$W$ is of type $B_n$ or $D_n$, $\o=1$ and $\cw=W$ with  $\cl|_{\cs}$ constant equal to $1$, or

$W$ is of type $D_n$, $n=2r\ge6$, $\o\in\Om''_W$, $\o^2=1$ and $\cw$ is of affine type $B_r$, with
the values of $\cl|_{\cs}$ being $1,2,2,\do,2$.

\subhead 5.5\endsubhead
We now list the various $(\cw,\cl)$ which are associated to various
$W,\o,W_J\in\cc_\o(W)$.

$\cw=\{1\}$ or $\cw$ is an irreducible affine Weyl groups of type $E_6,E_7,E_8$ or $D_m,m\ge4$ with $\cl|_{\cs}$ constant equal to $1$;

$\cw$ of affine type $A_{n-1},n\ge3$ with $\cl|_{\cs}$ constant in $\ZZ_{>0}$;

$\cw$ of affine type $B_m,m\ge3$ with the values of $\cl|_{\cs}$ being
$1,1,\do,1,1$ or $1,2,2,\do 2$;

$\cw$ of affine type $C_m,m\ge1$ with the values of $\cl|_{\cs}$ being
$t,1,1,\do,1,1,s$ with $(t,s)\in\ZZ_{>0}\T\ZZ_{>0}$ such that $t-s\in\{0,1,-1\}$;

$\cw$ of affine type $C_m,m\ge2$ with the values of $\cl|_{\cs}$ being
$t,2,2,\do,2,2,s$ with $(t,s)\in\ZZ_{>0}\T\ZZ_{>0}$ such that $t-s\in\{1,-1\}$;

$\cw$ of affine type $G_2$ with the values of $\cl|_{\cs}$ being $1,1,1$ or
$3,3,1$;

$\cw$ of affine type $F_4$ with the values of $\cl|_{\cs}$ being $1,1,1,1,1$ or
$2,2,2,1,1$.
\nl
(Here affine of type $C_1$ is taken to be the same as affine of type $A_1$.)

\subhead 5.6\endsubhead
Let $\o\in\Om_W$, $J\subsetneqq S$ be such that $W_J\in\cc_\o(W)$.
Let $E\m c_E$ be as in 4.1 with $\cw=\cw_J$.
For $E,\tE$ in $\Irr(\bar\cw_J)$ we write $E\le\tE$ if $E=\tE$ or $c_E>c_{\tE}$. This
is a partial order on $\Irr(\bar\cw_J)$.
For $E,\tE$ in $\Irr(\bar\cw_J)$ we write $E\approx\tE$ if $c_E=c_{\tE}$. This
is an equivalence relation on $\Irr(\bar\cw_J)$.

Let $\qq$ be an indeterminate. Let $\un\ct_J=\QQ(\qq)\ot\ct_J$.
Now the obvious $\cw_J$ action on $\ct_J$ induces a linear action of $\cw_J$ on
$\un\ct_J$.

For $E,\tE$ in $\Irr(\bar\cw_J)$ we define
$$\Om'_{E,\tE}=\sha(\bar\cw)\i\sum_{w\in\bar\cw_J}\tr(w\i,E)\tr(w,\tE)
\det(\qq-w,\un\ct_J)\i\qq^{-c_E-c_{\tE}}\in\QQ(\qq).$$
The following result can be deduced from Lemma 2.1 in \cite{GM} (where it is attributed
to the author).

\proclaim{Proposition 5.7} The system of equations
$$\sum_{E',\tE'\text{ in }\bar\cw_J;E'\le E,\tE'\le\tE,E'\si \tE'}
P_{E',E}\L'_{E',\tE'}P_{\tE',\tE}=\Om'_{E,\tE},$$
$$P_{E,E}=1 \text{for all }E,$$
$$P_{E',E}=0 \text{ if }E'\approx E,E\ne E',$$
$$P_{E',E}=0 \text{ if }E'\not\le E,$$
$$\L'_{E',E}=0 \text{ unless }E'\approx E$$
with unknowns
$$P_{E',E}\in\QQ(\qq),\L'_{E',E}\in\QQ(\qq),(E',E \text{ in }\Irr(\bar\cw_J))$$
has a unique solution. 
\endproclaim

\head 6. The set $\cg_\o(G)$ and the bijection $\cg_\o(G)@>\si>>\cc_\o(W)$ \endhead
\subhead 6.1\endsubhead
Let $le(G)$ be the set of subgroups of $G$ (see 0.1)
which are Levi subgroups of some parabolic subgroup of $G$.
Let $ls^0(G)$ be the subset of $ls(G)$ (see 0.1) consisting
of the pairs $(\fc,\fL)$ each of which is a unipotent block by itself.
Let $\un ls(G)$ be the set of unipotent blocks of $G$. 
Note that $ls^0(G)$ can be identified with a subset $\un ls^0(G)$ of $\un ls(G)$.
Let $Z_G$ be the group of components of the centre of $G$ and let
$Z_G^*=\Hom(Z_G,\CC^*)$.
We have a partition $ls(G)=\sqc_{\c\in Z_G^*}ls_\c(G)$ where $ls_\c(G)$
consists of all $(\fc,\fL)\in ls(G)$ such that the natural action of $Z_G$ on 
$\fL$ is through $\c$. For any $\c\in Z_G^*$, $ls_\c(C)$ is a union of
unipotent blocks. Hence we have a partition 
$$\un{ls}(G)=\sqc_{\c\in Z_G^*}\un{ls}_\c(G)$$ 
where $\un{ls}_\c(G)$ is the set of unipotent blocks contained in $ls_\c(G)$.
In \cite{L84} it is shown that for any $\c\in Z_G^*$ 

(a) the intersection $ls^0(G)\cap ls_\c(G)$ consists of at most one element. 
\nl
Let $\cg'(G)$ be the set consisting of $G$-conjugacy classes of triples 
$(L,\fc_1,\fL_1)$ with $L\in le(G)$, $(\fc_1,\fL_1)\in ls^0(L)$. 
In \cite{L84} a bijection 

(b) $\un{ls}(G)\lra\cg'(G)$
\nl
is established. 
To  $(L,\fc_1,\fL_1)\in\cg'(G)$ we have associated in \cite{L84, 4.4} a perverse
sheaf $\ph_!K$ (up to shift) 
on $G$ whose cohomology sheaves restricted to unipotent classes
are direct sums of local systems in a fixed unipotent block. (This defines the
unipotent blocks and the bijection (b).)
If $\c\in Z_G^*$ and 
$\b\in\un{ls}_\c(G)$, then the corresponding $(L,\fc_1,\fL_1)\in\cg'(G)$ is such 
that $(\fc_1,\fL_1)\in ls_{\c_1}(L)$ for some $\c_1\in Z_L^*$ which is 
uniquely determined by $\c$. (Note that $\c$ is the image of $\c_1$ under the
injective homomorphism $Z_L^*@>>>Z_G^*$ induced by the obvious (surjective)
homomorphism $Z_G@>>>Z_L$.) Using (a) for $L$ we see that
$(\fc_1,\fL_1)$ is unique if it exists. Thus (b) gives rise to a bijection

(c) $\un{ls}_\c(G)\lra\cg_\c(G)$
\nl
where $\cg_\c(G)$ is the set of all $L\in le(G)$ (up to $G$-conjugacy)
such that there exists $(\fc_1,\fL_1)\in ls^0(L)\cap ls_{\c_1}(L)$
where $\c_1\in Z_L^*$ maps to $\c$ under $Z_L^*@>>>Z_G^*$.
We set $\cg(G)=\sqc_{\c\in Z_G^*}\cg_\c(G)$. Then (c) gives rise to a
bijection

(d) $\un{ls}(G)\lra\cg(G)$.
\nl
In \cite{L84} it is shown
that the set of objects in a fixed unipotent block $\b$ is in bijection
with the set of irreducible representations of the normalizer $\fW_\b$ of $L$
modulo $L$ (with $L$ corresponding to $\b$) and that $\fW_\b$ is naturally a 
Weyl group.

\subhead 6.2 \endsubhead
In the remainder of this section we assume that $G$ in 0.1 is almost simple,
simply connected and that $W,S$ in 3.1 is the affine Weyl group associated with a
simple algebraic group $H$ over $\CC$ of type dual to that of $G$. In particular,
$\ct$ (see 3.1) can be identified with the group of one parameter subgroups
of a maximal torus $T$ of $H$.
For $h\in T$, the connected centralizer of $h$ in $H$ has Weyl group equal to
a $W$-conjugate of $W_J$ for some $J\subsetneqq S$. This gives a correspondence
$h---J$ between $T$ and $\{J;J\subsetneqq S\}$.

We can find an isomorphism $\io:Z_G^*@>\si>>\Om_W$
such that the following holds: if $G=Spin_N(\CC)$, with $N\ge5$ odd or
$N\ge10$ even, the subset $\Om''_W$ of $\Om_W$ (see 3.1) corresponds to the set
of characters of $Z_G$ which do not factor through $SO_N(\CC)$. We shall
identify $Z_G^*=\Om_W$ via $\io$.

The following result appears in \cite{L79}, \cite{L09}.
\proclaim{Theorem 6.3} Define $un(G)@>>>\Irr(\bW)$ by $\fc\m E$
where $E$ is attached to $(\fc,\CC)$ under the Springer correspondence. This map
defines a bijection $un(G)@>\si>>\bS(\bW)$ (notation of 3.2).
\endproclaim

\subhead 6.4\endsubhead
Let $\o\in\Om_W$. 
Let $L\in\cg_\o(G)$, let $(L,\fc_1,\fL_1)$ be the corresponding object of $\cg'_G$
(see 6.1) and let $\ph_!K$ be the associated complex of sheaves on $G$
(see 6.1); note that $\ph_!K[m]$ is a semisimple perverse sheaf for some $m$.
Let $K_1$ be a simple perverse sheaf on $G$ such that $K_1$ is a direct summand
of $\ph_!K[m]$. Then $K_1$ is a character sheaf on $G$. Let $C$ be the semisimple
conjugacy class of $H$ attached to $K_1$ by the known classification of
character sheaves. Note that $C$ is independent of the choice of $K_1$.
There is a unique subset $J\subsetneqq S$ such that for any $h\in T\cap C$
we have $h---J$ (see 6.2).

\proclaim{Theorem 6.5} For $L\in\cg_\o(G)$ we define $J\subsetneqq S$ as in 6.4. We have $W_J\in\cc_\o(W)$ and $L\m W_J$ is a bijection
$$\cg_\o(G)@>\si>>\cc_\o(W).\tag a$$
\endproclaim
Note that the Theorem gives a parametrization of the set of unipotent blocks of $G$
which is purely in terms of $W$ and is thus independent of the geometry of $G$.
The proof is given in the remainder of this section.

\subhead 6.6\endsubhead
Let $\o\in\Om_W$. We now describe explicitly (and independently of
6.4) a bijection
$$\cg_\o(G)@>>>\cc_\o(W).\tag a$$
The set $\cg_\o(G)$ is computed in \cite{L84}.
Assume first that $G=SL_n(\CC)$, $n\ge2$.
Note that $Z_G^*=\Om_W$ is a cyclic group of order $n$.
Let $k=ord(\o)$ (a divisor of $n$). Now $\cg_\o(G)$ 
consists of a single $L\in le(G)$ (up to conjugacy)
such that $L_{der}\cong SL_k(\CC)\T\do\T SL_k(\CC)$ ($n/k$ copies). 
The bijection (a) is the obvious bijection between sets with one element.

Assume that $G$ is of type $E_6$. Then $\Om_W$ is cyclic of order $3$.
If $\o=1$ then $\cg_\o(G)$ consists of a single object, a maximal torus $L$. The
bijection (a) is the obvious bijection between sets with one element. 
If $\o\ne1$ then $\cg_\o(G)$ consists of $L=G$ and of
$L\in le(G)$ such that $L_{der}=SL_3(\CC)\T SL_3(\CC)$ (up to conjugacy).
Recall that $\cc_\o(W)=\{\{1\},W_J\}$ where $W_J$ is of type $D_4$.
We define (a) by $G\m W_J$, $L\m\{1\}$ where $L\ne G$.

Assume that $G$ is of type $E_7$. Then $\Om_W$ is cyclic of order $2$.
If $\o=1$ then $\cg_\o(G)$ consists of a single object, a
maximal torus $L$. The bijection (a) is the obvious bijection between sets with
one element. If $\o\ne1$ then $\cg_\o(G)$ consists of $L=G$ and of
$L\in le(G)$ (up to conjugacy) such that
$L_{der}=SL_2(\CC)\T SL_2(\CC)\T SL_2(\CC)$ and $L$ is contained in an
$L'\in le(G)$ with $L'_{der}=SL_6(\CC)$ but $L$ is not contained in an
$L''\in le(G)$ with $L''_{der}=SL_7(\CC)$.
Recall that $\cc_\o(W)=\{\{1\},W_J\}$ where $W_J$ is of type $E_6$.
We define (a) by $G\m W_J$, $L\m\{1\}$ where $L\ne G$.

Assume that $G$ is of type $E_8,F_4$ or $G_2$. 
We have $\o=1$ and $\cg_\o(G)$ consists of two objects,
an $L$ which is a maximal torus and $G$.
Recall that $\cc_\o(W)=\{\{1\},W_J\}$ where $W_J$ is
of the non-affine type $E_8,F_4$ or $G_2$ (respectively).
We define (a) by $G\m W_J$, $L\m\{1\}$ where $L\ne G$.

\subhead 6.7\endsubhead
In this subsection we assume that $G$ is $Sp_{2n}(\CC), (n\ge3)$,
$Spin_{2n+1}(\CC),(n\ge2)$ or $Spin_{2n}(\CC), (n\ge4)$ so that
$W$ is of affine type $B_n (n\ge3),C_n (n\ge2)$ or $D_n (n\ge4)$ (respectively).

From \cite{L84} we see that $\cg_\o(G)$ consists of all
$L\in le(G)$ (up to conjugacy) such that

$L_{der}\cong Sp_{\d}$, for various $(\d,r)\in {}'\un\cc_\o(W)$ (for $W$ of type $B$),

$L_{der}\cong Spin_{\d}$, for various $(\d,r)\in{}'\un\cc_\o(W)$ (for $W$ of type $C$
or $D$, $\o\in\Om'_W$),

$L_{der}\cong Spin_{\d}\T SL_2(\CC)^r$, for various $(\d,r)\in{}'\un\cc_\o(W)$
(for $W$ of type $C$ or $D$, $\o\in\Om''_W$).

Then $L\m(\d,r)$ defines a bijection $\cg_\o(G)@>\si>>{}'\un\cc_\o(W)$. Composing
this with the inverses of the bijections ${}'\cc_\o(W)@>\si>>{}'\un\cc_\o(W)$
(3.10, 3.11) and $\cc_\o(W)@>\si>>{}'\cc_\o(W)$ (3.8, 3.9) we obtain the bijection
6.6(a) in our case.

This completes the definition of the bijection 6.6(a) in all cases.
It can be verified that this associates to $L\in \cg_\o(G)$ the same $W_J$
as that defined in 6.5. Hence the map 6.5(a) is well defined and it is a bijection
(the same as 6.6(a)).

\subhead 6.8\endsubhead
Let $\o\in\Om_W,\b\in\un{ls}_\o(G)$  with corresponding $L\in\cg_\o(G)$.
Under 6.5(a), $L$ corresponds to $W_J\in\cc_\o(W)$. Let
$(\cw,\cl,\cs,\bar\cw)=(\cw_J,\cl_J,\cs_J,\bar\cw_J)$ be as in \S5.

Let $(L,\fc_1,\fL_1)\in\cg'(G)$ be a triple corresponding to $L$ as in
6.1.
Let $\fc_{max}$ be the unipotent class of $G$ induced by $\fc_1$; let
$\fc_{min}$ be the unipotent class of $G$ that contains $\fc_1$.

For any $\fc\in un(G)$ we denote by $b_\fc$ the dimension of the variety of Borel
subgroups of $G$ that contain a fixed element of $\fc$. A case by case
verification gives the following two results.

\proclaim{Theorem 6.9} We have

(a) $b_{\fc_{max}}=\aa[{}^\o W_J]$ (notation of 2.1).

(b) $b_{\fc_{min}}-b_{\fc_{max}}=\nu(\cw,\cl)$ (notation of 4.3).
\endproclaim

\proclaim{Theorem 6.10}There exists a group isomorphism $\fW_\b@>\si>>\bar\cw$ well
defined up to composition with an inner automorphism of
$\bar\cw$ given by the action of an element in $\Om_\cw$ (see 3.1).
It carries the set of simple reflections of $\fW_\b$ into the image of $\cs$
under $\cw@>>>\bar\cw$.
\endproclaim

\proclaim{Conjecture 6.11} Assume that $\cw\ne\{1\}$. Let
$(\fc,\fL),(\fc',\fL')$ be in $\b$ and let $E,E'$ be in $\Irr(\fW_\b)=\Irr(\bar\cw)$
(this equality follows from 6.10). Assume that the generalized Springer
correspondence \cite{L84} associates $E$ to $(\fc,\fL)$ and $E'$ to $(\fc',\fL')$.
Then

(a) $b_\fc-b_{\fc_{max}}=c_E$ where $c_E$ is defined as in 4.1 in terms of
$(\cw,\cl)$;

(b) we have $E\si E'$ if and only if $\fc=\fc'$.
\endproclaim
This holds in the case where $G$ is of exceptional type. (In the case
where $G$ is of type $E_8,F_4$ or $G_2$ this follows from \cite{L20}.) This can
be also verified in the cases where $(\cw,\cl)$ is as in the examples in 4.2.

\subhead 6.12\endsubhead
For any $i=(\fc,\fL),i'=(\fc',\fL')$ in $\fb$ let $\Om_{i,i'}\in\QQ(\qq)$ be as
in \cite{L86, 24.7} and let $\Pi_{i',i}\in\QQ(\qq)$ be as in \cite{L86, 24.8}.
Here $\qq$ is an indeterminate. Let $E,E'$ in $\Irr(\bar\cw)$ be corresponding to
$i,i'$ (respectively)
under the generalized Springer correspondence. From the definitions
we have $\Om_{i,i'}=f(\qq)\Om'_{E,E'}$ where $f(\qq)\in\QQ(\qq)-\{0\}$ is
independent of $i,i'$. Assuming that 6.11(a) holds we see from  5.7
and \cite{L86, 24.8} that 
$$\Pi_{i',i}=P_{E',E}\tag a$$
for any $i,i'$ as above. In particular this holds in the case where $G$ is of
exceptional type. Note that $\Pi_{i',i}$ measures the stalks of the intersection
cohomology sheaf on the closure of $\fc$ with coefficients in $\fL$ while
$P_{E',E}$ is determined purely in terms of $W$.

\head 7. Cells in the weighted affine Weyl group $\cw_J$\endhead
\subhead 7.1\endsubhead
Let
$$\z:Cell(W,||)@>\si>>un(G)\tag a$$
be the bijection defined in \cite{L89}.

\subhead 7.2\endsubhead
We now fix $\o\in\Om_W$ and $\b\in\un{ls}_\o(G)$ with
corresponding $L\in\cg_\o(G)$. Under  6.5(a), $L$ corresponds to
$J\subsetneqq S$ such that $W_J\in\cc_\o(W)$.
Let $\cw_J,\cs_J,\cl_J$ be as in \S5. We assume that $\sha(\cs_J)\ge2$.
Let $un_\b(G)$ be the set of all $\fc\in un(G)$ such that $(\fc,\fL)\in\b$
for some $\fL$.
We define a map

(a) $\fA:Cell(\cw_J,\cl_J)@>>>un(G)$
\nl
asssuming a conjecture in \cite{L02, \S25}.
Let $\boc\in Cell(\cw_J,\cl_J)$.
Let $E_0$ be the unique irreducible special representation of $W_J$ such that
$z(E_0)=r(op)$ (notation of 2.1 with $\fW$ replaced by $W_J$) and let
$\boc_0\in Cell(W_J,||)$ be such that $E_0$ belongs to $\boc_0$. According to
\cite{L03, Conj.25.3} there is a well defined $\ti\boc\in Cell(W,||)$
which contains $yx$ for any
$y\in\boc_0,x\in\boc'$ (the product $yx$ is taken in $W$). (We use the fact that,
by results of \cite{L84a}, $W_J$ satisfies the assumptions of \cite{L03, 25.2}.)
 We set $\fA(\boc)=\z(\ti\boc)$, with $\z$ as in 7.1(a).

\proclaim{Conjecture 7.3}$\fA$ is injective with image equal to $un_\b(G)$.
Hence $\fA$ defines a bijection $Cell(\cw_J,\cl_J)@>\si>>un_\b(G)$.
\endproclaim
This is a generalization of 7.1(a).
Note that \cite{L03, 25.2} holds if $J=\emp$ in which case it states that
any $\boc\in Cell(\cw_J,\cl_J)$  is contained in a two-sided cell of $(W,||)$;
this can be deduced from \cite{L03, 10.14}. Using this one can define $\fA$
unconditionally when $W$ is of exceptional type and verify the conjecture
in that case.

\subhead 7.4\endsubhead
Let $\aa:\cw_J@>>>\NN$ be the $a$-function (see \cite{L03}) of the weighted
affine Weyl group $(\cw_J,\cl_J)$. Let $\boc\in Cell(\cw_J,\cl_J)$ and let
$\fc=\fA_J(\boc)\in un_\b(G)$. We expect that 

(a) for any $w\in\boc$ we have $\aa(w)=b_\fc-b_{\fc_{max}}$.

\enddocument